\documentclass[a4paper]{article}
\usepackage[latin1]{inputenc} 
\usepackage[T1]{fontenc} 
\usepackage{RR,RRthemes}
\usepackage{hyperref}

\usepackage[dvips]{epsfig,graphicx}
\usepackage{verbatim}
\usepackage{array}
\usepackage{color}
\usepackage{colortbl}
\usepackage{fancybox}
\usepackage{longtable}
\usepackage{epsfig}
\usepackage{amsmath,amsfonts}
\usepackage{amssymb}

\usepackage{latexsym}
\usepackage{graphicx}
\usepackage{rotating}
\usepackage{longtable,dcolumn,multirow}

\newcommand{\norm}[1]{\Vert {#1}\Vert }
\newcommand{\neuc}[1]{\left \Vert {#1} \right \Vert_2 }
\newcommand{\nfro}[1]{\Vert {#1}\Vert_F }
\newcommand{\dab}{(\Delta A,\Delta b)}
\newcommand{\dA}{\Delta A}
\newcommand{\db}{\Delta b}
\newcommand{\tx}{\tilde{x}}
\newcommand{\mcg}{{\cal {M}}_{g'}}
\newcommand{\finproof}{\begin{flushright} $\Box$ \end{flushright}}
\newcommand{\zc}{0_{n,1}}
\newcommand{\zr}{0_{1,n}}

\newtheorem{propos}{Proposition}
\newtheorem{theo}{Theorem}
\newtheorem{coro}{Corollary}

\newtheorem{Remark}{Remark}
\newtheorem{Algo}{Algorithm}

\RRdate{D\'ecembre 2010}

\RRauthor{
Marc Baboulin\thanks{Universit\'e Paris-Sud and INRIA, France
   ({\tt marc.baboulin@inria.fr}).}
  \and
Serge Gratton\thanks{ENSEEIHT and CERFACS, France
   ({\tt serge.gratton@enseeiht.fr}).}
}
\authorhead{Baboulin \& Gratton}
\RRtitle{Une contribution au conditionnement du probl\`eme de moindres carr\'es totaux}
\RRetitle{A contribution to the conditioning of the total least squares problem}
\titlehead{Conditioning of the total least squares}
\RRnote{Also appeared as LAPACK Working Note 236}
\RRresume{Nous d\'erivons des formules exactes pour le conditionnement d'une fonction lin\'eaire de la solution d'un probl\`eme de moindres carr\'es totaux. Etant donn\'e un syst\`eme lin\'eaire sur-d\'etermin\'e $Ax=b$, nous montrons que ce conditionnement peut \^etre calcul\'e en utilisant les valeurs singuli\`eres et les vecteurs singuliers \`a droite de $[A,b]$ et $A$.
Nous proposons aussi une borne sup\'erieure qui n\'ecessite le calcul de la plus grande et de la plus petite valeur singuli\`ere de $[A,b]$
ainsi que de la plus petite valeur singuli\`ere de $A$. Dans des exemples num\'eriques, nous comparons ces conditionnements et les 
bornes d'erreur directe correspondantes avec les erreurs estimatives donn\'ees dans ~\cite{vava91}.
}
\RRabstract{
We derive closed formulas for the condition number of a linear function of the total least squares solution.
Given an over determined linear system $Ax=b$, we show that this condition number can be
computed using the singular values and the right singular vectors of $[A,b]$ and $A$.
We also provide an upper bound that requires the computation of the largest and the smallest singular value of $[A,b]$
and the smallest singular value of $A$.
In numerical examples, we compare these values and the resulting forward error bounds
with the error estimates given in~\cite{vava91}.}
\RRmotcle{moindres carr\'es totaux, conditionnement, perturbations normwise,
mod\`ele d'erreur dans les variables}
\RRkeyword{total least squares, condition number, normwise perturbations, errors-in-variables model}
\RRprojet{Grand-Large}  
\RRdomaine{Informatique/Analyse num\'erique} 
\RCSaclay 

\begin{document}
\RRNo{7488}
\makeRR   

\bibliographystyle{plain}

\section{Introduction}
Given a matrix $A\in \Bbb{R}^{m \times n}~(m>n)$ and an observation vector $b \in \Bbb{R}^{m}$, the standard
over determined linear least squares (LS) problem consists in finding a vector $x \in \Bbb{R}^{n}$ such that $Ax$ is the best approximation of $b$.
Such a problem can be formulated using what is referred to as the linear statistical model
$$
b = Ax +\epsilon,
~A \in \mathbb{R}^{m\times n},
~b \in \mathbb{R}^m,
~{\rm rank}(A)=n,
$$
where $\epsilon$ is a vector of random errors having expected value
$E(\epsilon)=0$ and variance-covariance $V(\epsilon)=\sigma^2I$.

In the linear statistical model, random errors affect exclusively the observation vector $b$ while $A$ is considered as known exactly. However it is often more realistic to consider that measurement errors might also affect $A$.
This case is treated by the statistical model referred to as Errors-In-Variables model (see e.g~\cite[p. 230]{vava91} and ~\cite[p. 176]{BJORCK}),
where we have the relation
$$(A+E)x=b+\epsilon.$$
In general it is assumed in this model that the rows of $[E, \epsilon]$ are independently and identically distributed with common zero mean vector and common covariance matrix.
The corresponding linear algebra problem, discussed originally in~\cite{GVL.80}, is
called the Total Least Squares (TLS) problem and can be expressed as:
\begin{equation}\label{tls}
\min_{E,\epsilon}\nfro{(E,\epsilon)},~~(A+E)x=b+\epsilon,
\end{equation}
where $\nfro{\cdot}$ denotes the Frobenius matrix norm.
As mentioned in~\cite[p. 238]{vava91}, the TLS method enables us to obtain a more
accurate solution when entries of $A$ are perturbed under certain conditions.

In error analysis, condition numbers are considered as fundamental tools since they
measure the effect on the solution of small changes in the data.
In particular the conditioning of the least squares problem was extensively studied in the numerical linear algebra literature (see e.g~\cite{BJORCK,CHA.IPS.95,COHI.99.1,CDW.07,ELDEN,GR.96,HIGHAM,IPSEN,KEN.LAU.98,STEWART}).
The more general case of the conditioning of a linear function of an LS solution was studied
in~\cite{ABG.07} and~\cite{BG.09} when perturbations on data are measured respectively normwise and componentwise
(note that the componentwise and normwise condition numbers for LS problems were also treated in~\cite{CDW.07} but without the generalization to a linear function of the solution).
Moreover we can find in~\cite{BDGL.09} algorithms
using the software libraries LAPACK~\cite{LAPACK} and ScaLAPACK~\cite{SCALAPACK}
as well as physical applications.

The notion of {\bf Total Least Squares} was initially defined in the seminal paper~\cite{GVL.80} that was the first to propose a numerically stable algorithm. Then various aspects of the TLS problem were developed in the comprehensive book~\cite{vava91} including a large survey of theoretical bases, computational methods and applications but also sensitivity analysis with for instance upper bounds for the TLS perturbation.
The so-called {\bf Scaled Total Least Squares} (STLS) problem
($\min_{E,\epsilon}\nfro{(E,\epsilon)},~~(A+E)x\gamma =\gamma b+\epsilon$,
for a given scaling parameter $\gamma$)
was formulated in~\cite{CORE} in which were addressed the difficulties coming from non existence of TLS solution.
In a recent paper~\cite{WEI.09}, we can find sharp estimates of the normwise, mixed and componentwise
condition numbers of the Scaled Total Least Squares (STLS) problem.

Here we are concerned with the TLS problem, which is a special case of the STLS problem,
and we will consider perturbations on data $(A,b)$ that are measured normwise using a product norm.
Contrary to~\cite{WEI.09}, we will consider the general case of the conditioning of $L^Tx$,
linear function of the TLS solution for which we will derive an exact formula.
The common situations correspond to the special cases where $L$ is the identity matrix
(condition number of the TLS solution) or a canonical vector
(condition number of one solution component).
The conditioning of a nonlinear function of a TLS solution can also be obtained
by replacing, in the condition number expression, the quantity $L^T$ by the Jacobian matrix at the solution.

We notice that the normwise condition number expression proposed in~\cite{WEI.09} is based on the evaluation of the norm of a matrix expressed as a Kronecker product resulting in large matrices which may be, as pointed out by the authors, impractical to compute, especially for large size problems.
We propose here a computable expression for the resulting condition number (exact formula and upper bound) using data that could be already available from the TLS solution process, namely by-products of the SVD decomposition of $A$ and $[A,b]$.
We also make use of the adjoint operator which enables us to work on a space of lower dimension
and we propose a practical algorithm based on the power method.

\section{Definitions and notations}
\subsection{The total least squares problem}

Let $A\in \Bbb{R}^{m \times n}$ and $b \in \Bbb{R}^{m}$, with $m>n$.
Following~\cite{vava91}, we consider the two singular value decompositions
of $A$, and $[A,b]$ : $A=U'\Sigma' V^{'T}$ and
$[A,b]=U\Sigma V^ T$. We also set $\Sigma=\mbox{diag}(\sigma_1,\dots,
\sigma_{n+1})$,
$\Sigma'=\mbox{diag}(\sigma'_1,\dots,\sigma'_n)$,
where the singular values are in nonincreasing order,
and define $\lambda_i=\sigma_i^2$, and $\lambda_i'=\sigma_i^{'2}$.
From~\cite[p. 178]{BJORCK}, we have the interlacing property
\begin{equation}\label{interlace}
\sigma_1 \geq \sigma'_1 \geq \sigma_2 \geq \dots \geq \sigma_n \geq \sigma'_n \geq \sigma_{n+1}.
\end{equation}
We consider the total least squares problem expressed in Equation~(\ref{tls})
%
and we assume in this text that the {\em genericity } condition $\sigma_n' > \sigma_{n+1}$ holds
(for more information about the "nongeneric" problem see e.g~\cite{vava91,CORE}).
From~\cite[Theorems~2.6 and~2.7]{vava91}, it follows that
the TLS solution $x$ exists, is unique, and satisfies
\begin{equation}\label{soltls}
x=\left (A^TA-\lambda_{n+1}I_n\right )^{-1}A^Tb.
\end{equation}
In addition, $\begin{bmatrix} x \\ -1 \end{bmatrix}$ is an eigenvector of
$[A,b]^T[A,b]$ associated with the simple eigenvalue $\lambda_{n+1}$,
i.e $\sigma_n' > \sigma_{n+1}$ guarantees that $\lambda_{n+1}$ is not a
semi-simple eigenvalue of $[A,b]^T[A,b]$. As for linear least squares problems,
we define the total least squares residual
$r=b-Ax$, which enables us to write
\begin{equation}\label{expvp}
\lambda_{n+1}=\frac{1}{1+x^Tx}\begin{bmatrix} x^T, & -1 \end{bmatrix}
\begin{bmatrix} A^TA & A^Tb \\ b^TA & b^Tb \end{bmatrix}
\begin{bmatrix} x \\ -1 \end{bmatrix}=\frac{r^Tr}{1+x^Tx}.
\end{equation}
As mentioned in~\cite[p. 35]{vava91}, the TLS solution is obtained by scaling the last right singular vector $v_{n+1}$ of
$[A,b]$ until its last.
component is $-1$ and, if $v_{i,n+1}$ denotes the $i$th component of $v_{n+1}$, we have
\begin{equation}\label{exp:basicsol}
x=-\frac{1}{v_{n+1,n+1}}[v_{1,n+1},\dots,v_{n,n+1}]^T.
\end{equation}
The TLS method involves an SVD computation and the
computational cost is higher than that of a classical LS problem
(about $2mn^2+12n^3$ as mentioned in~\cite[p. 598]{GOLUB},
to be compared with the approximately $2mn^2$ flops
required for LS solved via Householder QR factorization). However, there exist faster methods
referred to as "partial SVD" (PSVD) that calculate only the last right singular vector or a basis
of the right singular subspace associated with the smallest singular values of $[A,b]$
(see~\cite[p. 97]{vava91}).
\subsection{Condition number of the TLS problem}
To measure the perturbations on data $A$ and $b$, we consider
the product norm defined on $\mathbb{R} ^{m \times n} \times \mathbb{R} ^m$ by  $\norm{(A,b)}_F=\sqrt{\norm{A}_{\rm{F}}^2+\neuc{b}^2}$
and we take the Euclidean norm $\neuc{x}$ for the solution space $\mathbb{R} ^n$.
In the following, the $n \times n$ identity matrix is denoted by $I_n$.

Let $L$ be a given $n \times k$ matrix, with $k \leq n$. We suppose here that $L$ is not
perturbed numerically and we  consider  the mapping
$$
\begin{array}{c c c c}
g\ : &
\mathbb{R} ^{m \times n} \times \mathbb{R} ^m &  \longrightarrow &  \mathbb{R} ^k
\\
 & ( A,b) & \longmapsto &  g(A,b)=L^{T}x=L^T(A^TA-\lambda_{n+1}I_n)^{-1}A^Tb,\\
 \end{array}
$$
Since $\lambda_{n+1}$ is simple, $g$ is a Fr\'echet-differentiable
function of $A$ and $b$, and the genericity assumption ensures that
the matrix $(A^TA-\lambda_{n+1}I_n)^{-1}$ is also Fr\'echet-differentiable
in a neighborhood of $(A,b)$. As a result, $g$ is Fr\'echet-differentiable
in a neighborhood of $(A,b)$.

The approach that we follow here is based on the work by~\cite{geur82,rice.66} where the mathematical difficulty of a problem is measured by the norm of the Fr\'echet derivative of the problem solution expressed as a function of data. This measure is an attainable bound at first order, and may therefore be approximate when large perturbations are considered.

Using the definition given in~\cite{geur82,rice.66},
we can express the condition number
of $L^Tx$, linear function of the TLS solution as
\begin{equation}
\label{exp:defcn}
K(L,A,b)
=\max_{\dab \ne 0}
\frac{\neuc{g'(A,b).\dab}}{\nfro{\dab}}.
\end{equation}

$K(L,A,b)$ is sometimes called the {\it absolute} condition number of $L^Tx$ as opposed to the {\it relative} condition number of $L^Tx$ and defined, when $L^Tx$ is nonzero by
\begin{equation}
\label{exp:defcnrel}
K^{(rel)}(L,A,b)=K(L,A,b) \nfro{(A,b)} / \neuc{L^Tx}.
\end{equation}

In the remainder, the quantity $K(L,A,b)$ will be simply referred to as the TLS condition number,
even though the proper conditioning of the TLS solution
corresponds to the special case when $L$ is the identity matrix.\\
In the expression $g'(A,b).\dab$, the "." operator denotes that we apply the linear function $g'(A,b)$ to the
variable $\dab$. We will use this notation throughout this paper to designate the image of a vector or a matrix by a linear function.
\begin{Remark}
{\rm
The case where $g(A,b)=h(x)$, with $h$ being a
differentiable nonlinear function mapping $\Bbb{R}^n$ to $\Bbb{R}^k$
is also covered because we have
$$g'(A,b).\dab=h'(x).(x'(A,b).\dab),$$
and $L^T$ would correspond to the Jacobian matrix $h'(x)$.
The nonlinear function $h$ can be for instance the Euclidean norm of part of the solution
(e.g in the computation of Fourier coefficients when we are interested in the quantity of signal in a given frequency band).
}
\end{Remark}
\section{Explicit formula for the TLS condition number}
\subsection{Fr\'echet derivative}
In this section, we compute the Fr\'echet d\'erivative of $g$ under
the genericity assumption, which enables us to obtain an explicit formula for the TLS condition number
in Proposition~\ref{propcn}.
\begin{propos}
\label{prop:der}
Under the genericity assumption,  $g$ is Fr\'echet differentiable
in a neighborhood of $(A,b)$. Setting $B_{\lambda}=A^TA-\lambda_{n+1}I_n$,
the Fr\'echet derivative
of $g$ at $(A,b)$ is expressed by
\begin{equation}
\label{funct}
\begin{array}{r  c l}
g'(A,b)\ :
\mathbb{R} ^{m \times n} \times \mathbb{R} ^m   \longrightarrow &  \mathbb{R} ^k
\\
 ( \dA,\db)  \longmapsto &
L^TB_{\lambda}^{-1}\left ( A^T +\frac{2x r^T}{1+x^Tx}\right )\left ( \db -\dA x\right ) +  \\
    &
L^TB_{\lambda}^{-1} \dA^T r.\\
 \end{array}
 \end{equation}
\end{propos}
{
\underline{Proof:}
The result is obtained from the chain rule.
Since $\lambda_{n+1}$, expressed in Equation~(\ref{expvp}), is a simple eigenvalue
of $[A,b]^T[A,b]$ with corresponding unit
eigenvector $\frac{1}{\sqrt{1+x^Tx}}\begin{bmatrix} x^T & -1 \end{bmatrix}^T$,
$\lambda_{n+1}$ is differentiable in a neighborhood of $(A,b)$ and then we have
\begin{eqnarray*}
\lambda'_{n+1}(A,b).\dab & = &
\frac{1}{1+x^Tx}\begin{bmatrix} x^T & -1 \end{bmatrix}
\begin{bmatrix} \dA^TA+A^T\dA & \dA^Tb+A^T\db \\ b^T\dA +\db^TA & \db^Tb+b^T\db \end{bmatrix}
\begin{bmatrix} x \\ -1 \end{bmatrix} \\
& = & \frac{2}{1+x^Tx}  \left (  x^T \dA^TAx -  x^T \dA^Tb -  x^T A^T\db + b^T\db \right ) \\
& = & \frac{2}{1+x^Tx}  \left ( -x^T \dA^T r +  (b^T - x^T A^T)\db \right ) \\
& = & \frac{2}{1+x^Tx}  \left ( -r^T \dA x +  r^T\db \right ),
\end{eqnarray*}
yielding
\begin{equation}
\label{eq:l}
\lambda'_{n+1}(A,b).\dab = \frac{2r^T(\db - \dA x)}{1+x^Tx}.
\end{equation}

Applying the chain rule to $B_{\lambda}^{-1}$, we obtain
\begin{eqnarray*}
(B_{\lambda}^{-1})'(A,b).\dab & = & -B_{\lambda}^{-1}\left ( \dA^TA+A^T\dA - \lambda'_{n+1}(A,b).\dab I_n \right ) B_{\lambda}^{-1} \nonumber\\
& = &  -B_{\lambda}^{-1}\left (  \dA^TA+A^T\dA - \frac{2r^T(\db - \dA x)}{1+x^Tx} I_n \right ) B_{\lambda}^{-1}.
\end{eqnarray*}

The chain rule now applied to $g(A,b)$ leads to
\begin{eqnarray*}
g'(A,b).\dab  & = & -L^T B_{\lambda}^{-1}\left ( \dA^TA+A^T\dA - \lambda'_{n+1}(A,b).\dab I_n \right ) B_{\lambda}^{-1}A^Tb+  L^T B_{\lambda}^{-1} \dA^T b + L^T B_{\lambda}^{-1} A ^T \db\\
& = &-L^T B_{\lambda}^{-1}\left ( \dA^TA+A^T\dA - \lambda'_{n+1}(A,b).\dab I_n \right ) x + L^T B_{\lambda}^{-1} \left (\dA^T b +  A ^T \db \right)\\
& = &
L^T B_{\lambda}^{-1}\left ( A^T +\frac{2x r^T}{1+x^Tx}\right )\left ( \db -\dA x\right )+L^T B_{\lambda}^{-1} \dA^T r,
\end{eqnarray*}
which gives the result.
\finproof
}

We now introduce the vec operation that stacks all the columns of a matrix into
a long vector:  for $A=[a_1,\dots,a_n]\in \Bbb{R}^{m \times n}$, vec$(A)=[a_1^T,\dots,a_n^T]^T \in \Bbb{R}^{mn \times 1}$.
Let $P \in \Bbb{R}^{mn  \times mn} $ denote the permutation matrix that represents
the matrix transpose by $vec(B^T)= P vec(B)$.  We remind also that
$vec(AXB)= (B^T \otimes A) vec(X)$, where $\otimes$ denotes the Kronecker product of two matrices
~\cite[p. 21]{KRONECKER}.

Let us now express the matrix representing $g'(A,b)$, denoted by $\mcg$.
Since $g'(A,b).(\dA,\db) \in \mathbb{R} ^k$, we have
$g'(A,b).(\dA,\db)=vec(g'(A,b).(\dA,\db))$ and setting in addition $D_\lambda=L^TB_{\lambda}^{-1}\left ( A^T +\frac{2x r^T}{1+x^Tx}\right )\in \Bbb{R}^{ k\times m}$, we obtain from~(\ref{funct})
\begin{eqnarray*}
g'(A,b).(\dA,\db)  & = &\mbox{vec}\left ( D_\lambda \left ( \db -\dA x\right )+L^TB_{\lambda}^{-1} \dA^T r\right ) \\
& = &\left(-x^T \otimes D_\lambda \right) \mbox{vec}(\dA)+\left( r^T \otimes (L^TB_{\lambda}^{-1}) \right) \mbox{vec}(\dA^T)+D_\lambda \db \\
& = &
\begin{bmatrix}
-x^T \otimes D_\lambda + \left(r^T \otimes (L^TB_{\lambda}^{-1})\right)P  ,&
D_\lambda
\end{bmatrix}
\begin{bmatrix}
\mbox{vec}(\dA) \\
\db
\end{bmatrix}.
\end{eqnarray*}
Then we get
$$\mcg =
\begin{bmatrix}
-x^T \otimes D_\lambda + \left(r^T \otimes (L^TB_{\lambda}^{-1})\right)P  ,&
D_\lambda
\end{bmatrix}
\in  \Bbb{R}^{k \times (nm+m)}.
$$
But we have $\norm{(\dA , \db)}_F =\neuc{
 \begin{bmatrix}
\mbox{vec}(\dA )  \\
\db
\end{bmatrix}
} $ and then, from Proposition~\ref{prop:der}
and using the definition of $K(L,A,b)$ given in Expression~(\ref{exp:defcn}),
we get the following proposition that expresses the TLS condition number in terms of the norm of a matrix.
\begin{propos}
\label{propcn}
The condition number of $g(A,b)$ is given by
$$K(L,A,b) = \neuc{\mcg},$$
where
\begin{eqnarray*}
\mcg & = &
\begin{bmatrix}
-x^T \otimes D_\lambda + \left(r^T \otimes (L^TB_{\lambda}^{-1})\right)P  ,&
D_\lambda
\end{bmatrix}
\in  \Bbb{R}^{k \times (nm+m)}.
\end{eqnarray*}
\end{propos}
\subsection{Adjoint operator and algorithm}\label{sec:adj}
Computing $K(L,A,b)$ reduces to computing the spectral norm
of the $k \times (nm+m)$ matrix $\mcg$. For large values of
$n$ or $m$, it is not possible to build explicitly the generally
dense matrix $\mcg$. Iterative techniques based on the power method~\cite[p. 289]{HIGHAM} or on the
Lanczos  method~\cite{GOLUB} are better suited. These algorithms
involve however the computation of the product of $\mcg^T$ by a
vector $y \in \Bbb{R}^k$. We describe now how to perform this operation.

Using successively the fact that $B_{\lambda}^{-T}=B_{\lambda}^{-1}$, $(A \otimes B)^T=A^T \otimes B^T$,
$vec(AXB)= (B^T \otimes A) vec(X)$ and $P^T=P^{-1}$ we have
\begin{eqnarray*}
\mcg^Ty & = &
\begin{bmatrix}
-x \otimes D_\lambda^T + P^T\left(r \otimes (B_{\lambda}^{-T}L)\right)     \\
D_\lambda^T
\end{bmatrix}
y \\
 & = &
\begin{bmatrix}
-(x \otimes D_\lambda^T) {\mbox{vec}}(y) + P^T\left(r \otimes (B_{\lambda}^{-1}L)\right) {\mbox{vec}}(y) \\
D_\lambda^T y
\end{bmatrix} \\
 & = &
\begin{bmatrix}
P^{-1} \left(P{\mbox{vec}} \left(-D_\lambda^T y x^T\right) + {\mbox{vec}}\left(B_{\lambda}^{-1}L y r^T\right) \right)\\
D_\lambda^T y
\end{bmatrix} \\
 & = &
\begin{bmatrix}
P^{-1} \left({\mbox{vec}} \left((-D_\lambda^T y x^T)^T\right) + {\mbox{vec}}\left(B_{\lambda}^{-1}L y r^T\right) \right)\\
D_\lambda^T y
\end{bmatrix} \\
 & = &
\begin{bmatrix}
P^{-1} {\mbox{vec}} \left(-x y^T D_\lambda + B_{\lambda}^{-1}L y r^T\right)\\
D_\lambda^T y
\end{bmatrix},
\end{eqnarray*}
and since for any matrix $B$ we have $P^{-1} vec(B)= vec(B^T)$, we get
\begin{equation}\label{exp:adj}
\mcg^Ty =
\begin{bmatrix}
{\mbox{vec}} \left (-D_\lambda^T y x^T + r y^T L^T B_{\lambda}^{-1} \right ) \\
D_\lambda^T y
\end{bmatrix}.
\end{equation}

This leads us to the following proposition.
\begin{propos}
\label{prop:adj}
The adjoint operator of $g'(A,b)$ using the scalar products\\
${\rm trace }(A_1^TA_2)+b_1^Tb_2$ and $y_1^Ty_2$ respectively on
$\mathbb{R} ^{m \times n} \times \mathbb{R} ^m$ and $\mathbb{R} ^k$ is
\begin{equation}
\label{functt}
\begin{array}{r  c l}
g^{'*}(A,b)\ :
 \mathbb{R} ^k \longrightarrow & \mathbb{R} ^{m \times n} \times \mathbb{R} ^m \\
y \longmapsto & \begin{pmatrix}
-D_\lambda^T y x^T + r y^T L^T B_{\lambda}^{-1},D_\lambda^T y
\end{pmatrix}
 \end{array}
 \end{equation}
In addition,  if $k=1$ we have
\begin{equation}
\label{kap1}
K(L,A,b)= \sqrt{\nfro{-D_\lambda^T  x^T + r  L^T B_{\lambda}^{-1}}^2+\neuc{D_\lambda}^2}
\end{equation}
\end{propos}

{
\underline{Proof:}
Let us denote by  $<(A_1,b_1),(A_2,b_2)>$ the scalar product ${\rm trace }(A_1^TA_2)+b_1^Tb_2$ on $\mathbb{R} ^{m \times n} \times \mathbb{R} ^m$.
We have  for any $y \in \Bbb{R}^k$,
\begin{eqnarray*}
y^T(g'(A,b).(\dA, \db))
& = &
y^T
\mcg \begin{bmatrix}
\mbox{vec}(\dA) \\
\db
\end{bmatrix} \\
& = &
(\mcg ^Ty)^T
\begin{bmatrix}
\mbox{vec}(\dA) \\
\db
\end{bmatrix}  \\
& = &
{\mbox{vec}} \left (-D_\lambda^T y x^T + r y^T L^T B_{\lambda}^{-1} \right )^T\mbox{vec}(\dA) + ( D_\lambda^T y)^T\db.
\end{eqnarray*}
Using now the fact that, for matrices $A_1$ and $A_2$ of identical sizes,\\
${\mbox{vec}} (A_1)^T {\mbox{vec}} (A_2)={\rm trace} (A_1^TA_2)$,
we get
\begin{eqnarray*}
y^T(g'(A,b).(\dA, \db)) & = &
{\rm trace} \left((-D_\lambda^T y x^T + r y^T L^T B_{\lambda}^{-1})^T \dA \right)
+ (D_\lambda^T y)^T\db\\
& = &
< \begin{pmatrix}
-D_\lambda^T y x^T + r y^T L^T B_{\lambda}^{-1},D_\lambda^T y
\end{pmatrix}
,(\dA, \db) >\\
& = &
<g^{'*}(A,b).y ,(\dA, \db) >,
\end{eqnarray*}
which concludes the first part of the proof.

For the second part, we
use
$$K(L,A,b)=\neuc{\mcg}=\neuc{\mcg^T}=\max_{y\ne 0}
\frac{
\neuc{
\begin{bmatrix}
{\mbox{vec}} \left (-D_\lambda^T y x^T + r y^T L^T B_{\lambda}^{-1} \right ) \\
D_\lambda^T y
\end{bmatrix}
}
}
{\Vert y\Vert_2}
$$
Since $k=1$, we have $y \in \Bbb{R}$, and
$K(L,A,b)=\neuc{
\begin{bmatrix}
{\mbox{vec}} \left (-D_\lambda^T  x^T + r  L^T B_{\lambda}^{-1} \right ) \\
{\mbox{vec}} ( D_\lambda^T ),
\end{bmatrix}
}
$
and the result follows from the relation ${\mbox{vec}} (A_1)^T {\mbox{vec}} (A_1)={\rm trace} A_1^TA_1=\nfro{A_1}^2.$
\finproof
}
\begin{Remark}
{\rm
The special case $k=1$ recovers the situation where we compute the conditioning
of the $i$th solution component.
In that case $L$ is the $i$th canonical vector of $\mathbb{R}^n$ and, in
Equation~(\ref{kap1}),
$L^T B_{\lambda}^{-1}$ is the $i$th row of $B_{\lambda}^{-1}$ and
$D_\lambda$ is the $i$th row of $B_{\lambda}^{-1}\left ( A^T +\frac{2x r^T}{1+x^Tx}\right )$.
}
\end{Remark}

Using~(\ref{funct}) and~(\ref{functt}), we can now write in Algorithm 1
the iteration of the power method (~\cite[p. 289]{HIGHAM}) to compute the TLS condition number $K(L,A,b)$.

\setlength{\fboxsep}{0.45cm}
\setlength{\fboxrule}{0.02cm}
\begin{center}
\fbox{
\begin{minipage}{0.8\textwidth}
{\bf Algorithm 1} : Condition number of TLS problem
\begin{description}
\item[~] Select initial vector $y \in \mathbb{R} ^k$
\item[~] {\bf for p=1,2,...}
\item[~] \hspace{0.5cm} $(A_p,b_p)=
\begin{pmatrix}-D_\lambda^T y x^T + r y^T L^T B_{\lambda}^{-1},D_\lambda^T y
\end{pmatrix}$
\item[~] \hspace{0.5cm} $\nu = \nfro{(A_p,b_p)} $
\item[~] \hspace{0.5cm} $(A_p,b_p)  \leftarrow ( \frac{1}{\nu}\cdot A_p,\frac{1}{\nu}\cdot b_p) $
\item[~] \hspace{0.5cm} $y= L^TB_{\lambda}^{-1}\left ( A^T +\frac{2x r^T}{1+x^Tx}\right )\left
( b_p - A_p x\right ) + L^TB_{\lambda}^{-1} A_p^T r $
\item[~] {\bf end}
\item[~] $K(L,A,b)=\sqrt{\nu}$
\end{description}
\end{minipage}
}
\end{center}
The quantity $\nu$ computed by Algorithm 1 is the largest eigenvalue of ${\mcg}{\mcg^T}$.
Since $K(L,A,b)=\neuc{\mcg}$ then the condition number $K(L,A,b)$ is also the largest singular value of ${\mcg}$ i.e $\sqrt{\nu}$.
As mentioned in~\cite[p. 331]{GOLUB}, the algorithm will converge if the initial $y$ has a component in the direction of the corresponding dominant eigenvector of ${\mcg}{\mcg^T}$. When there is an estimate of this dominant eigenvector, the initial $y$ can be set to this estimate but in many implementations, $y$ is initialized as a random vector.
The algorithm is terminated by a "sufficiently" large number of iterations or by
evaluating the difference between two successive values of $\nu$ and comparing it to a tolerance given by the user.

\subsection{Closed formula}
Using the adjoint formulas obtained in Section~\ref{sec:adj}, we now get a closed formula for the total least squares conditioning.
\begin{theo}
\label{theo:cond}
We consider the total least squares problem and assume that the genericity assumption holds.
Setting $B_{\lambda} = A^TA-\lambda_{n+1}I_n$,
then the condition number of $L^Tx$, linear function of the TLS solution, is expressed by
$$
K(L,A,b)=\neuc{C}^{\frac{1}{2}},$$
where $C$ is the $k \times k$ symmetric matrix
$$
C=(1+\neuc{x}^2) L^TB_{\lambda}^{-1}
\left (
A^TA+
\lambda_{n+1}(I_n-\frac{2 x x^T}{1+\neuc{x}^2})
\right)
B_{\lambda}^{-1}L
.$$
\end{theo}
{
\underline{Proof:}
We have
$K(L,A,b)^2=\neuc{\mcg^T}^2=\max_{\neuc{y}=1}
\neuc{\mcg^Ty}^2$. If $y$ is a unit vector in $\Bbb{R}^{k}$, then using Equation~(\ref{exp:adj}) we obtain
\begin{eqnarray*}
\neuc{\mcg^Ty}^2 & = &
\neuc{{\rm vec}\left (-D_\lambda^T y x^T + r y^T L^T B_{\lambda}^{-1} \right ) }^2 + \neuc{D_{\lambda}^Ty}^2\\
&=& \nfro{-D_\lambda^T y x^T + r y^T L^T B_{\lambda}^{-1}  }^2 + \neuc{D_{\lambda}^Ty}^2\\
&=&   \nfro{D_\lambda^T y x^T}^2+\nfro{ r y^T L^T B_{\lambda}^{-1}}^2-
2\, { \rm trace}( x y^T D_\lambda  r y^T L^T B_{\lambda}^{-1}) + \neuc{D_{\lambda}^Ty}^2.
\end{eqnarray*}
For all vectors $u$ and $v$, we have $\nfro{uv^T}=\neuc{u}\neuc{v}$. Moreover we have
$$
{ \rm trace}\left( (x y^T D_\lambda  r) (y^T L^T B_{\lambda}^{-1})\right)
= { \rm trace}\left((y^T L^T B_{\lambda}^{-1}) (x y^T D_\lambda  r) \right)
= y^T L^T B_{\lambda}^{-1} x  r^T D_{\lambda}^Ty.$$
Thus
\begin{eqnarray*}
\neuc{\mcg^Ty}^2 & = &
\neuc{x}^2\neuc{D_\lambda^T y}^2 + \neuc{r}^2 \neuc{B_{\lambda}^{-1}Ly}^2 -2\,
y^T L^T B_{\lambda}^{-1} x r^T D_\lambda^T y + \neuc{D_{\lambda}^Ty}^2\\
& = &
(1+x^Tx)y^T D_{\lambda} D_{\lambda}^T y+\neuc{r}^2 y^T L^T B_{\lambda}^{-2} L y -2\,
y^T L^T B_{\lambda}^{-1} x r^T D_\lambda^T y\\
& = &
y^T \left((1+x^Tx)D_{\lambda} D_{\lambda}^T + \neuc{r}^2 L^T B_{\lambda}^{-2} L -2\,
L^T B_{\lambda}^{-1} x r^T D_\lambda^T\right) y,
\end{eqnarray*}
i.e
$\neuc{\mcg^T}^2=\neuc{C}$ with
\begin{equation}\label{exp1:C}
C=(1+x^Tx)D_{\lambda} D_{\lambda}^T + \neuc{r}^2 L^T B_{\lambda}^{-2} L -2\,
L^T B_{\lambda}^{-1} x r^T D_\lambda^T.
\end{equation}
Replacing $D_{\lambda}$ by
$L^TB_{\lambda}^{-1}\left( A^T +\frac{2x r^T}{1+x^Tx}\right)$, Equation~(\ref{exp1:C}) simplifies to
\begin{equation}\label{exp2:C}
C=L^TB_{\lambda}^{-1}\left(
(1+x^Tx)A^TA+\neuc{r}^2I_n+2A^Trx^T
\right)B_{\lambda}^{-1}L.
\end{equation}
But $A^Trx^T=A^T(b-Ax)x^T=A^Tbx^T-A^TAxx^T$ and, since from Equation~(\ref{soltls}) we have $A^Tb=B_{\lambda}x$,
we get
$A^Trx^T=B_{\lambda}xx^T-A^TAxx^T=(A^TA-\lambda_{n+1}I_n)xx^T-A^TAxx^T=-\lambda_{n+1}xx^T.$
From Equation~(\ref{expvp}) we also have $\neuc{r}^2=\lambda_{n+1}(1+x^Tx)$ and thus
Equation~(\ref{exp2:C}) becomes
\begin{eqnarray*}
C & = & L^TB_{\lambda}^{-1}\left(
(1+x^Tx)A^TA+\lambda_{n+1}(1+x^Tx)I_n-2\lambda_{n+1}xx^T
\right)B_{\lambda}^{-1}L\\
& = &
(1+\neuc{x}^2) L^TB_{\lambda}^{-1}
\left (
A^TA+
\lambda_{n+1}(I_n-\frac{2 x x^T}{1+\neuc{x}^2})
\right)
B_{\lambda}^{-1}L.
\end{eqnarray*}
\finproof
}
\section{TLS condition number and SVD}
\subsection{Closed formula and upper bound}

Computing $K(L,A,b)$ using Theorem~\ref{theo:cond} requires the explicit formation of the normal equations matrix $A^TA$ which is a source of rounding errors and also generates an extra computational cost of about $mn^2$ flops. In practice the TLS solution is obtained by Equation~(\ref{exp:basicsol}) and involves an SVD computation.
In the following theorem, we propose a formula for $K(L,A,b)$ that can be computed with quantities that may be already available from the solution process.
In the following $\zc$ (resp. $\zr$) denotes the zero column (resp. row) vector of length $n$.
\begin{theo}
\label{theo:cntls}
Let $V$ and $V'$ be the matrices whose columns are the right singular vectors of respectively $[A,b]$ and $A$ associated with the singular values
$(\sigma_1,\dots,\sigma_{n+1})$ and $(\sigma'_1,\dots,\sigma'_n)$. Then the condition number of $L^Tx$, linear function of the TLS solution is expressed by
$$
K(L,A,b)=(1+\neuc{x}^2)^{\frac{1}{2}}
\neuc{
L^T V' D'
\begin{bmatrix}
V^{'T}, &
\zc
\end{bmatrix}
V
\begin{bmatrix}
D, &
\zc
\end{bmatrix}^T
}
,~{\rm where}~$$
$D'=\mbox{diag}\left((\sigma_1^{'2}-\sigma_{n+1}^2)^{-1},\dots,(\sigma_n^{'2}-\sigma_{n+1}^2)^{-1}\right)
~{\rm and}~
D=\mbox{diag}\left((\sigma_1^2+\sigma_{n+1}^2)^{\frac{1}{2}}
,\dots,(\sigma_n^2+\sigma_{n+1}^2)^{\frac{1}{2}}\right).
$
When $L$ is the identity matrix, then the condition number reduces to
$$
K(L,A,b)=(1+\neuc{x}^2)^{\frac{1}{2}}
\neuc{
D'
\begin{bmatrix}
V^{'T}, &
\zc
\end{bmatrix}
V
\begin{bmatrix}
D, &
\zc
\end{bmatrix}^T}.$$
\end{theo}

{
\underline{Proof:}
From $[A,b]=U\Sigma V^T$, we have
$[A,b]^T[A,b]=V \Sigma^2 V^T=\sum_{i=1}^{n+1} \sigma_i^2 v_i v_i^T$ and
\begin{eqnarray*}
[A,b]^T[A,b]+\lambda_{n+1}I_{n+1} & = & \sum_{i=1}^{n+1} \sigma_i^2 v_i v_i^T
+ \lambda_{n+1} \sum_{i=1}^{n+1} v_i v_i^T\\
 & = & \sum_{i=1}^{n+1}(\sigma_i^2 + \lambda_{n+1}) v_i v_i^T\\
 & = & \sum_{i=1}^{n}(\sigma_i^2 + \sigma_{n+1}^2) v_i v_i^T + 2 \lambda_{n+1} v_{n+1} v_{n+1}^T,
\end{eqnarray*}
leading to
\begin{equation}\label{exp:demcntls1}
[A,b]^T[A,b]+\lambda_{n+1}I_{n+1}-2 \lambda_{n+1} v_{n+1} v_{n+1}^T =
\sum_{i=1}^{n}(\sigma_i^2 + \sigma_{n+1}^2) v_i v_i^T
\end{equation}
From Equation~(\ref{exp:basicsol}), we have $v_{n+1}=-v_{n+1,n+1}\begin{bmatrix} x \\ -1 \end{bmatrix}$
and, since $v_{n+1}$ is a unit vector, $v_{n+1,n+1}^2=\frac{1}{1+\neuc{x}^2}$.
Then Equation~(\ref{exp:demcntls1}) can be expressed in matrix notation as
\begin{equation}\label{exp:demcntls2}
\begin{bmatrix} A^TA & A^Tb \\ b^TA & b^Tb \end{bmatrix}+
\lambda_{n+1} \begin{bmatrix} I_n & \zc \\ \zr & 1 \end{bmatrix}
-\frac{2 \lambda_{n+1}}{1+\neuc{x}^2} \begin{bmatrix} xx^T & -x \\ -x^T & 1 \end{bmatrix}
=\sum_{i=1}^{n}(\sigma_i^2 + \sigma_{n+1}^2) v_i v_i^T
\end{equation}
The quantity $A^TA+\lambda_{n+1}(I_n-\frac{2 x x^T}{1+\neuc{x}^2})$
corresponds to the left-hand side of Equation~(\ref{exp:demcntls2}) in which the last row and the last column have
been removed. Thus it can also be written
$$
A^TA+\lambda_{n+1}(I_n-\frac{2 x x^T}{1+\neuc{x}^2})
=
\begin{bmatrix} I_n, & \zc \end{bmatrix}
\left( \sum_{i=1}^{n}(\sigma_i^2 + \sigma_{n+1}^2) v_i v_i^T \right)
\begin{bmatrix} I_n \\ \zr \end{bmatrix},
$$
and the matrix $C$ from Theorem~\ref{theo:cond} can be expressed
\begin{equation}\label{exp:demcntls3}
C=(1+\neuc{x}^2)L^T
\begin{bmatrix} B_{\lambda}^{-1}, & \zc \end{bmatrix}
\left( \sum_{i=1}^{n}(\sigma_i^2 + \sigma_{n+1}^2) v_i v_i^T \right)
\begin{bmatrix} B_{\lambda}^{-1} \\ \zr \end{bmatrix}
L.
\end{equation}
Moreover from $A=U'\Sigma' V^{'T}$, we have
$A^TA=V' \Sigma^{'2} V^{'T}=\sum_{i=1}^{n} \sigma_i^{'2} v'_i v_i^{'T}$ and
\begin{eqnarray*}
 B_{\lambda} & = & A^TA-\lambda_{n+1}I_n\\
 & = & \sum_{i=1}^{n}\sigma_i^{'2} v'_i v_i^{'T} - \sigma_{n+1}^2 \sum_{i=1}^{n} v'_i v_i^{'T}\\
 & = & \sum_{i=1}^{n} (\sigma_i^{'2}-\sigma_{n+1}^2 ) v'_i v_i^{'T}\\
 & = & V' D'^{-1} V^{'T}.
\end{eqnarray*}
Hence $B_{\lambda}^{-1}=V^{'-T} D' V^{'-1}=V' D' V^{'T}$ and
$\begin{bmatrix} B_{\lambda}^{-1}, & \zc \end{bmatrix}=V' D' \begin{bmatrix}  V^{'T}, & \zc \end{bmatrix}$.\\
We also have $\sum_{i=1}^{n}(\sigma_i^2 + \sigma_{n+1}^2) v_i v_i^T
=V \begin{bmatrix} D \\ \zr \end{bmatrix}
\begin{bmatrix} D, & \zc \end{bmatrix} V^T.
$\\
Then, by replacing in Equation~(\ref{exp:demcntls3}), we obtain
$C=(1+\neuc{x}^2) \widetilde{V} \widetilde{V}^T$ with
$\widetilde{V}=L^T V' D'
\begin{bmatrix}
V^{'T}, &
\zc
\end{bmatrix}
V
\begin{bmatrix}
D, &
\zc
\end{bmatrix}^T
$.
As a result, using Theorem~\ref{theo:cond},
$$K(L,A,b)^2=
\neuc{C}=(1+\neuc{x}^2) \neuc{\widetilde{V} \widetilde{V}^T}=(1+\neuc{x}^2) \neuc{\widetilde{V}}^2.$$
When $L=I_n$, we use the fact that $V'$ is an orthogonal matrix and can be removed from the expression of
$\neuc{\widetilde{V}}^2$.
\finproof
In many applications, an upper bound would be sufficient to give an estimate of the conditioning of the TLS solution. The following corollary gives an upper bound for $K(L,A,b)$.
\begin{coro}
\label{coro:ubcntls}
The condition number of $L^Tx$, linear function of the TLS solution is bounded by
$$
{\bar{K}}(L,A,b)= (1+\neuc{x}^2)^{\frac{1}{2}} \neuc{L}
\frac{(\sigma_1^2+\sigma_{n+1}^2)^{\frac{1}{2}}}{(\sigma_n^{'2}-\sigma_{n+1}^2)}.
$$
\end{coro}
{
\underline{Proof:}
This result comes from the inequality
$\neuc{AB} \leq\neuc{A} \neuc{B}$, followed by
$\neuc{D'}=\max_{i} (\sigma_i^{'2}-\sigma_{n+1}^2)^{-1} =
(\sigma_n^{'2}-\sigma_{n+1}^2)^{-1}$ and
$\neuc{D}^2=\max_{i} (\sigma_i^2+\sigma_{n+1}^2)=
(\sigma_1^2+\sigma_{n+1}^2)$.
\finproof

\subsection{Numerical examples}
In the following examples we study the condition number of $x$ i.e $L$ is here the identity matrix $I_n$.
Then, to simplify the notations, we removed the variable $L$ from the expressions and the condition number of $x$
 will be denoted by $K(A,b)$ and its upper bound by ${\bar{K}}(A,b)$.
All the experiments were performed with MATLAB 7.6.0 using a machine precision $2.22\cdot10^{-16}$.
\subsubsection{First example}
In the first example we consider the TLS problem $Ax \approx b$ where $[A,b]$ is defined by
$$[A,b]=Y
\left(
\begin{array}{c}
D\\
0\\
\end{array}
\right)
Z^T \in \mathbb{R}^{m\times (n+1)}, Y=I_m-2yy^T, Z=I_{n+1}-2zz^T,$$
where $y \in \mathbb{R}^m$ and $z \in \mathbb{R}^{n+1}$ are random unit vectors,
$D=diag(n,n-1,\cdots,1,1-e_p)$ for a given parameter $e_p$.
The quantity $\sigma'_n-\sigma_{n+1}$ measures the distance of our problem to
nongenericity and, due to Equation~(\ref{interlace}), we have in exact arithmetic
$$\sigma'_n-\sigma_{n+1} \leq \sigma_n-\sigma_{n+1} = e_p.$$
Then by varying $e_p$, we can generate different TLS problems and by considering small values of $e_p$, it is possible to study the behavior of the TLS condition number in the context of close-to-nongeneric problems.
The TLS solution $x$ is computed using an SVD of $[A,b]$ and Equation~(\ref{exp:basicsol}).

In Table~\ref{tab:num2}, we compare the exact condition number $K(A,b)$ given in Theorem~\ref{theo:cntls},
the upper bound ${\bar{K}}(A,b)$ given in Corollary~\ref{coro:ubcntls},
and the upper bound obtained from~\cite[p. 212]{vava91}
and expressed by
$$
\kappa(A,b)= \frac{9\sigma_1 \neuc{x}}{\sigma_n-\sigma_{n+1}} \left(1+\frac{\neuc{b}}{\sigma'_n-\sigma_{n+1}} \right)
\frac{1}{\neuc{b}-\sigma_{n+1}}.
$$
We also report the condition number computed by Algorithm 1, denoted by $K_p(A,b)$,
and the corresponding number of power iterations
(the algorithm terminates when the difference between two successive values is lower than $10^{-8}$).
When $\sigma'_n-\sigma_{n+1}$ decreases, the TLS problem becomes worse conditioned and there is a factor ${\cal O}(10)$
between the exact condition number $K(A,b)$ and its upper bound ${\bar{K}}(A,b)$.
We also observe that ${\bar{K}}(A,b)$ is an estimate of better order of magnitude than
$\kappa(A,b)$ and that, for small values of $\sigma'_n-\sigma_{n+1}$, $\kappa(A,b)$ is much less reliable.
$K_p(A,b)$ is always equal or very close to $K(A,b)$.
\begin{table}[hbtp!]
\centering
\caption{TLS conditioning for several values of $\sigma'_n-\sigma_{n+1}$.}
\vspace{0.4cm}
\begin{tabular}{|c|c|c|c|c|c|}
\hline
&&&&&\\
$\sigma'_n-\sigma_{n+1}$&$K(A,b)$&${\bar{K}}(A,b)$&$\kappa(A,b)$&$K_p(A,b)$&\#iter\\
&&&&&\\
\hline
$9.99976032\cdot10^{-1}$&$1.18\cdot10^{0}$&$2.36\cdot10^{1}$&$1.29\cdot10^{2}$&$1.18\cdot10^{0}$&$11$\\
\hline
$9.99952397\cdot10^{-5}$&$8.36\cdot10^{3}$&$1.18\cdot10^{5}$&$1.31\cdot10^{10}$&$8.36\cdot10^{3}$&$6$\\
\hline
$9.99952365\cdot10^{-9}$&$8.36\cdot10^{7}$&$1.18\cdot10^{9}$&$1.31\cdot10^{18}$&$8.36\cdot10^{7}$&$4$\\
\hline
$9.99644811\cdot10^{-13}$&$8.36\cdot10^{11}$&$1.18\cdot10^{13}$&$1.31\cdot10^{26}$&$8.32\cdot10^{11}$&$5$\\
\hline
\end{tabular}
\label{tab:num2}
\end{table}
\subsubsection{Second example}
Let us now consider the following example from~\cite[p. 42]{vava91} also used in ~\cite{WEI.09} where
$$A=\left(
\begin{array}{cccc}
m-1&-1&\cdots&-1\\
-1&m-1&\cdots&-1\\
\cdot&&&\\
\cdot&&&\\
\cdot&&&\\
-1&-1&\cdots&m-1\\
-1&-1&\cdots&-1\\
-1&-1&\cdots&-1\\
\end{array}
\right) \in \mathbb{R}^{m\times (m-2)},
~b=\left(
\begin{array}{c}
-1\\
-1\\
\cdot\\
\cdot\\
\cdot\\
-1\\
m-1\\
-1\\
\end{array}
\right) \in \mathbb{R}^{m}.
$$
The exact solution of the TLS problem $Ax \approx b$ can be computed analytically~\cite[p. 42]{vava91} and is equal to $x=-(1,\cdots,1)^T$.
We consider a random perturbation $(\dA,\db)$ of small norm $\nfro{(\dA,\db)}=10^{-10}$ and we denote by $\tx$ the computed solution of the perturbed system $(A+\dA) x \approx b+\db$.

In Table~\ref{tab:num1}, we report for several values of $m$ the {\it relative} condition number as defined in~(\ref{exp:defcnrel})
and we compare the computed relative forward error $\frac{\neuc{\tx-x}}{\neuc{x}}$ with the forward error bounds that can be expected from the computation of $K^{(rel)}(A,b)$ and its upper bounds ${\bar{K}}^{(rel)}(A,b)$ and $\kappa^{(rel)}(A,b)$.
Since the condition number corresponds to the worst case in error amplification at first order, these quantities are, as observed in Table~\ref{tab:num1}, always larger than the computed forward error
(there is approximately a factor $10^2$ between those quantities).
We also observe that, in this example, ${\bar{K}}^{(rel)}(A,b)$ and $\kappa^{(rel)}(A,b)$ produce forward error estimates that are of same order of magnitude.
\begin{table}[hbtp!]
\centering
\caption{Forward error and upper bounds for a perturbed TLS problem.}
\vspace{0.4cm}
\begin{tabular}{|c|c|c|c|c|c|c|}
\hline
&&&&&\\
$m$&$K^{(rel)}(A,b)$&$\frac{\neuc{\tx-x}}{\neuc{x}}$&$K^{(rel)}(A,b) \frac{\nfro{(\dA,\db)}}{\nfro{(A,b)}}$ &${\bar{K}}^{(rel)}(A,b)\frac{\nfro{(\dA,\db)}}{\nfro{(A,b)}}$&$\kappa^{(rel)}(A,b)\frac{\nfro{(\dA,\db)}}{\nfro{(A,b)}}$\\
&&&&&\\
\hline
50&$5.05\cdot10^{1}$&$2.45\cdot10^{-13}$&$2.21\cdot10^{-11}$&$1.55\cdot10^{-10}$&$6.72\cdot10^{-10}$\\
\hline
100&$1.01\cdot10^{2}$&$1.08\cdot10^{-13}$&$1.55\cdot10^{-11}$&$1.54\cdot10^{-10}$&$4.26\cdot10^{-10}$\\
\hline
500&$5.01\cdot10^{2}$&$8.79\cdot10^{-14}$&$6.85\cdot10^{-12}$&$1.53\cdot10^{-10}$&$1.66\cdot10^{-10}$\\
\hline
1000&$1.00\cdot10^{3}$&$4.33\cdot10^{-14}$&$4.84\cdot10^{-12}$&$1.53\cdot10^{-10}$&$1.13\cdot10^{-10}$\\
\hline
\end{tabular}
\label{tab:num1}
\end{table}
\section{Conclusion}
We proposed sensitivity analysis tools for the total least squares problem when the genericity condition is satisfied.
We provided closed formulas for the condition number of a linear function of the TLS solution when the perturbations
of data are measured normwise.
We also described an algorithm based on an adjoint formula and we expressed this condition number and an upper bound
of it in terms of the SVDs of $[A,b]$ and $A$. We illustrated the use for these quantities in two numerical examples.

\bibliography{biblio}

\newpage
\tableofcontents

\end{document}